\documentclass[reqno]{amsart}
\usepackage{tikz-cd}
\usepackage{amsmath,amssymb,latexsym,graphicx}
\numberwithin{equation}{section}
\def\pf{{\bf Proof }}
\newtheorem{theorem}{Theorem}[section]
\newtheorem{lemma}[theorem]{Lemma}

\newtheorem{remark}[theorem]{Remark}

\def\qed{$\Box$}

\def\pf{{\bf Proof }}
\usepackage{enumerate}

\begin{document}

\title{Convex bodies with algebraic section volume functions}

%\title{On algebraically integrable bodies}

%    Only \author and \address are required; other information is
%    optional.  Remove any unused author tags.

%    author one information
\author[M. Agranovsky]{MARK AGRANOVSKY}
%\author{Mark Agranovsky}
\address{Bar-Ilan University and Holon Institute of Technology}
%\curraddr{}
\email{agranovs@math.biu.ac.il}
%\thanks{The author thanks Mikhail Zaidenberg for drawing the author's attention to the subject  and  stimulating discussions, and Victor %Vassiliev for useful correspondence.}

%    author two information
%\author{}

%\address{}
%\curraddr{}
%\email{}
%\thanks{}

%    The 2010 edition of the Mathematics Subject Classification is
%    the current definitive version.

%\date{}

\maketitle

\begin{abstract}
 The section volume function $A_K(\xi,t), \ \xi \in \mathbb R^n, \ t \in \mathbb R,$ of a body $K \subset \mathbb R^n$  evaluates the $(n-1)$-dimensional volume of the cross-section $K$ by the hyperplane $\{ x \cdot \xi=t \}.$  We are concerned with the question: can the shape of a body $K$ be detected from an algebraic type of its section function? We prove that among strictly convex bodies $K$ with $C^{\infty}$ boundaries, ellipsoids are completely described by the algebraic equation
 $qA_K^m+p=0,$ where $m \in \mathbb N$ and  $q=q(\xi), \ p=p(\xi,t)$ are polynomials.  The result is motivated by Arnold's problem on algebraically integrable domains (which, in turn, has its roots in Newton's Lemma about ovals), and generalizes known results (\cite{KMY}, \cite{A1}), \cite{A3}) on polynomially integrable domains.

\end{abstract}

%\footnote{MSC 2010: 44A12, 51M99; keywords: convex domains, $X$-tray transform, algebraic hypersurface, ellipsoid. }

\section{Introduction}\label{S:Intro}

\subsection{Formulation of the problem}
In geometric tomography, one of the important metric characteristics of a body $K \subset \mathbb R^n$ (a compact set with nonempty interior) is the
section volume function (or, for brevity, the {\it section function} \cite{Gardner}) :
$$A_K(\xi,t)= Vol_{n-1} \big[ K \cap X_{\xi,t} \big],$$
evaluating the $(n-1)$-dimensional volume of the cross-section the body $K$ by the hyperplane
$$X(\xi,t)=  \{x \in \mathbb R^n: x \cdot \xi=t \}$$
with the normal vector $ \ \xi \in S^{n-1}$, at the signed distance $t$ to the origin.  Here $ x \cdot \xi$ is the inner product, $S^{n-1}$ is the unit sphere in $\mathbb R^n.$

Equivalently, the section function can be  defined as the Radon transform
$$A_K(\xi,t)=(\mathcal R\chi_K)(\xi,t)$$
of the characteristic function $\chi_K$ of the body $K.$

The natural domain of definition of the function $A_K(\xi,t)$ consists of such pairs $(\xi,t)$ that the hyperplane $X(\xi,t)$ has a nonempty intersection with the interior of $K.$ For strictly convex domains, it is equivalent to  $A_K(\xi,t)>0.$

In \cite{A1} the notion of polynomially integrable bodies was introduced. The body $K$ is {\it polynomially integrable} if $A_K(\xi,t)$ is a polynomial in $t$ as long as $X(\xi,t) \cap K \neq \emptyset.$  It was proved in \cite{A1} that all polynomially integrable domains with $C^{\infty}$ boundaries are convex, there is no domains with these properties
if the dimension $n$ is even, and if $n$ is odd and the degree of the polynomial $t \to A_K(\xi,t)$ does not exceed $n-1$ then $\partial K$ is an ellipsoid.
Koldobsky, Merkurjev and Yaskin \cite{KMY} got rid of the restriction on the degree of the polynomial and therefore completely characterized polynomially integrable domains.

For the ellipsoids, the section function $A_K(\xi,t)$ is a polynomial in $t$ in odd dimensions, and  is the square root of a polynomial in even dimensions. Therefore,
in both cases, if $K$ is an ellipsoid then $A_K^2(\xi,t)$ is a polynomial in $t.$ It was asked in \cite{A3} whether the converse is true?
A partial answer was given in  \cite{AKRY} where it was proved that if $A_K^2(\xi,t)=P^2(\xi,t)Q(\xi,t),$ where $P, Q$ are polynomials in $t$ and $deg Q \leq 2$, then  $K$ is an ellipsoid. The conjecture in \cite{A3} is that, in any dimension, ellipsoids are characterized by the property that some power $A_K^m(\xi,t)$ is a polynomial in $t.$ This conjecture was confirmed in \cite{A3} for  $n=2$ under the assumption that the domain is bounded by a real algebraic curve. In this case $A_K$ is the $X$-ray transform of  the characteristic function $\chi_K.$ In this article, we extend the result of \cite{A3} to higher dimensions and confirm the conjecture in \cite{A3}, under an  assumption that the algebraic condition for $A_K(\xi,t)$ extends also to the variable $\xi.$  In fact, the algebraic condition in $\xi$ can be omitted if to assume from the very beginning that
$\partial K$ is a real algebraic hypersurface in $\mathbb R^n.$

\subsection{Main result}

The main result of this article is the following
\begin{theorem} \label{T:main} Let $K \subset \mathbb R^n$ be a strictly convex compact body with $C^{\infty}$ boundary $\partial K.$  Suppose that there exists $m \in \mathbb N,$ and nonzero polynomials polynomial $q(x)$ and $p(x,t)$ in $\mathbb R^n \times \mathbb R$ such that the algebraic equation
\begin{equation}\label{E:equation}
q(\xi)A_K^m(\xi,t)+p(\xi,t)=0
\end{equation}
is satisfied for all $\xi \in S^{n-1}, \ t \in \mathbb R$ such that $A_K(\xi,t) > 0.$
Then $\partial K$ is an ellipsoid.
\end{theorem}

\begin{remark} Since the function $A_K(\xi,t)$ is homogeneous of degree 0, the coefficients $p, q$ can be taken nonzero homogeneous polynomials of equal degrees.
\end{remark}
\begin{remark}
If $K=E$ is an ellipsoid, written in the standard form \
$E=\{ \sum\limits_{j=1}^n \frac{x_j^2}{a_j^2} \leq 1 \},$
then
$A_E(\xi,t)= const \ h_E^{-n}(\xi) (h_E^2(\xi)-t^2)^{\frac{n-1}{2}} \ , $
where  $h_E(\xi)= \big( \sum\limits_{j=1}^n a_j^2 \xi_j^2 \big)^{\frac{1}{2}} $ is the support function. It follows that for ellipsoids
equation \ref{E:equation} is fulfilled with  $m=2$ and $q(\xi)=h_E^{2n}(\xi), \ p(\xi,t)= (h_E^2(\xi)-t^2)^{n-1}.$
\end{remark}
\begin{remark}
As it will be seen from the proof of Theorem \ref{T:main}, the polynomiality of the coefficients $p, q$ with respect to the full variable $(\xi,t)$ is used
only for proving that the boundary $\partial K$ is  semi-algebraic. Therefore, the conclusion of Theorem \ref{T:main} remains true if $p(\xi,t)$ is assumed to be a polynomial only with respect to $t$ but the boundary $\partial K$ is a priori assumed semi-algebraic, i.e., $\partial K$ is a connected component of the zero set of a nonzero real polynomial. For such class of domains, Theorem \ref {T:main} states that $K$ is an ellipsoid if $A_K^m(\xi,t)$ is a polynomial in $t$ for some $m\in \mathbb N.$ This generalizes, for domains with algebraic boundaries, the results of \cite{KMY} and \cite{AKRY} on polynomially integrable domains.
\end{remark}

\subsection{Motivation and comments}
In 1987, V. Arnold asked (\cite{Ar},  1987-14, 1988-13, 1990-27) for higher-dimensional generalization of Lemma 28 of Newton about ovals
\cite{Principia, Wiki}. Newton proved, in connection with celestial mechanics, that the areas of the two parts  into which  a domain, bounded by a smooth oval in the plane, is divided by a straight line never depend
algebraically on the parameters of the secant line  (see \cite{ArKvant}).

Arnold introduced the notion of algebraically integrable domains in $\mathbb R^n$ which are those that the volumes of the portions cut off by a hyperplane (cutoff two-valued
function) algebraically depend on the parameters of the secant hyperplane. He conjectured that such domains with $C^{\infty}$ boundaries never exist if $n$ is even,  and all such domains in odd dimensions are ellipsoids. V. A. Vassiliev (see \cite{ArVas, Vas1, Vas2}) proved that there is no algebraically integrable $C^{\infty}$ domains in even-dimensional spaces. The "odd" part of Arnold's conjecture remains unanswered and only partial results are known by today \cite{Vas2, Vas3}.

Past decade, a series of results \cite{A1, A2, A3, A4, AKRY, B1, B2, KMY, YM} (see also the survey \cite{ABKVY}) has been obtained where the integrability problem is addressed the section function $A_K(\xi,t)$ rather than the  cutoff function. Notice that the condition of algebraicity addressed the section function is much weaker than that for the cutoff function and therefore the algebraic equation for the section function should be specified.
Polynomially integrable domains $K$ satisfy, by the definition, the equation $w+p=0,$  for $w=A_K(\xi,t),$ where $p=p(\xi,t)$ is a polynomial in $t$.
In this article we consider more general family of algebraic equations of order $m:$
$qw^m+p=0, \ w=A_K(\xi,t) >0.$
The result of this article is that the above two-terms non-linear algebraic equation completely describes ellipsoids in Euclidean spaces.

\section{The outline of the proof}

The key ingredient of the proof  is proving the upper bound \\
$\deg p(\xi, \cdot) \leq m(n-1)$ for the degree of the polynomial $p(\xi,t)$ in (\ref{E:equation}) with respect to the variable $t.$ The degree of a polynomial is determined by its growth at infinity. However, the behavior of $p(\xi,t)$ cannot be traced within the real space, because the information about $p(\xi,t)$ is available only in bounded interval of $t \in \mathbb R$, namely, for such $t$ that the hyperplane $X(\xi,t)$ meets $K$. This difficulty can be avoided by extending the section function for complex values of $t$ and studying the behavior of $p(\xi,t)$  when $t$ goes to infinity along a complex path.
The complexification of the whole construction is possible because the algebraic equation for $A_K(\xi,t)$  implies that
the boundary $\partial K$ is an algebraic hypersurface and hence extends, in a natural way,  to $\mathbb C^n$ as a complex algebraic hypersurface $\Gamma$.
Correspondingly, the section function $A_K(\xi,t)$ can be extended to complex values of $t.$

In more details, we construct, following  \cite{Vas1, Vas2}, the  analytic extension of the section function onto the complex Grassmanian
by integration of a certain differential volume $(n-1)$-form
on real $(n-1)$-chains in the complex space. The integrals coincide on homologically equivalent, relatively to $\Gamma,$ $(n-1)$-chains and hence define for fixed $\xi$ and $t,$ a function on the relative homology group. To construct a single-valued function of $\xi,t$ which analytically extends $A_K(\xi,t)$ in a containing $\infty$ domain of complex $t,$  we consider  a locally trivial homological vector bundle, which assigns to each $t \in \mathbb C$ a  fiber of the homology classes, on which the integrals of the above differential form are naturally defined. The restriction of the bundle to a smaller, avoiding branching points, unbounded simply-connected domain in the complex $t$-plane trivializes the above bundle  and the desired single-valued analytic extension of $A_K(\xi,t)$ is delivered by choosing a continuous section of the restricted bundle. Here we make use of the fact that all the objects are algebraic, and exploit the projective model by passing to homogeneous coordinates. It makes easier the analysis of the behavior of the analytic extension at infinity, which becomes an "infinitely distant"  hyperplane in the projective realization. This construction leads to proving the asymptotic $O(t^{n-1}), \ t \to \infty,$ for the analytic extension of $A_K(\xi,t),$ as $t \to \infty$ along a complex path, and yields the upper bound $deg_t p(\xi,t) \leq m(n-1),$  due to (\ref{E:equation}).

The final step of the proof is as follows (see \cite{A1}). We prove that at almost all points of $\partial K$  the polynomial $p(\xi,t)$ has two zeros at two parallel tangent hyperplanes $t=h_K(\xi), \ t=-h_K(-\xi),$  each zero has the multiplicity  $\frac{m(n-1)}{2}.$  Together with the estimate $deg_t p(\xi,t) \leq m(n-1)$ this leads to the conclusion that the above two multiple zeros  are the only roots of the polynomial $p(\xi,t).$  Bezout theorem gives
a simple expression of $p(\xi,t),$ and therefore, of the section function $A_K(\xi,t),$ in terms of the support function $h_K(\xi).$
To complete the proof, we use that $A_K=\mathcal R(\chi_K)$ belongs to the range of the Radon transform $\mathcal R$ and apply
the range conditions \cite{He} (the moment conditions with respect to $t$) for the Radon transform.  The combination of the first three moment conditions and the above explicit expression of $A_K(\xi,t)$ via the support function $h_K(\xi)$ results in a system of functional equations which implies that $h^2_K(\xi)$ is a quadratic polynomial. This is equivalent to $\partial K$ being an ellipsoid.

\section{Preparations}

\subsection{Algebraization and complexification}

From now on,  we fix  a compact connected strictly convex domain $K \subset \mathbb R^n$ with $C^{\infty}$ boundary $\partial K$ and with section function satisfying equation (\ref{E:equation}) in Theorem \ref{T:main}.

The next lemma shows that equation (\ref{E:equation}) implies that the boundary $\partial K$ is {\it semialgebraic}. By definition, it means that $\partial K$ is defined by finitely many polynomial equalities and inequalities. In our case, it means that $\partial K$ is a smooth connected component of the zero set of a real non-zero irreducible polynomial $Q.$

\begin{lemma} The boundary $\partial K$ is a semialgebraic non-singular hypersurface, i.e. there exists a nonzero irreducible polynomial $Q(x_1, \cdots, x_n)$
such that $\partial K$ is a smooth connected component of the zero set $Q^{-1}(0).$ Also $\partial K$ is infinitely smooth.
\end{lemma}

\pf
To show that $\partial K$ is semialgebraic, we follow arguments in \cite{Vas1, Vas2, ArVas, ArKvant} and refer to these sources for more details.
Let $\xi =\nu_a \in S^{n-1}$ be the unit normal vector at a point $a \in \partial K.$ Since $\partial K$ is strictly convex, the intersection of the tangent plane  $T_a(\partial K) =\{ (x-a) \cdot \nu_a=0 \}$ with the domain $K$ consists of the single point $a$   and hence $A_K(\nu_a, \nu_a \cdot a)=0.$  Then equation (\ref{E:equation}) yields $p(\nu_a, \nu_a \cdot a)=0.$
Therefore, the tangent bundle $(a,T_a(\partial K))_{ a \in \partial K} $ belongs to the zero set of a non-zero polynomial and hence is semialgebraic.
Then  $\partial K$ is semialgebraic, by projective duality.

By the condition, $\partial K \in C^{\infty}$ and hence $\partial K$ is a non-singular semialgebraic hypersurface in $\mathbb R^n$,  i.e.,
 $\partial K \subset Q^{-1}(0), \ Q$ is a nonzero irreducible polynomial, $\nabla Q(x) \neq 0, x \in \partial K, $ and $\partial K \in C^{\infty}.$

\qed

Everywhere in the sequel, $Q$ denotes a nonzero irreducible polynomial of $n$ variables with real coefficients such that
$$\partial K \subset \{x \in \mathbb R^n: Q(x)=0 \},$$
and $\partial K$ is a connected smooth component of the full zero set $Q^{-1}(0)$ in $\mathbb R^n.$

The polynomial $Q(x), x \in \mathbb R^n,$ extends, in a natural way, as a polynomial $Q(z)$ in the complex space $\mathbb C^n.$

Denote $\Gamma$ the  complex algebraic variety
$$\Gamma=\{ z \in \mathbb C^n: Q(z)=0\}.$$

In other words, $\Gamma$ is the complexification of the real algebraic variety $Q^{-1}(0) \subset \mathbb R^n.$
We will use the decomposition of the polynomial $Q$ into homogeneous polynomials
\begin{equation} \label{E:homo}
Q= Q_d+Q_{d-1}+ \cdots + Q_1+Q_0, \ d=dim \ Q.
\end{equation}
The zero variety
$$\Gamma_{\infty}=Q_d^{-1}(0)$$
of the leading homogenous term is {\it the asymptotic cone} of $\Gamma.$

The complex algebraic hypersurface $\Gamma$ can have singular points $z$ in $\mathbb C^n.$ They are points on $\Gamma$ which are  critical points of the polynomial $Q,$
that is $ z \in \Gamma$ with   $\frac{\partial Q}{\partial z_j}(z)=0, \ j=1,\cdots, n.$

We will  also be considering complex hyperplanes
$$X^{\mathbb C}(\xi,t)=\{ z \in \mathbb C^n: \xi_1 z_1 + \cdots + \xi_n z_n=t \} ,$$
where the parameters $\xi, t$ can be  complex numbers.

In the sequel, we will need the following
\begin{lemma}\label{L:Zariski}  There is a Zariski open subset $\Xi \subset S^{n-1}$  such that for any $\xi \in \Xi$  holds:
\begin{enumerate}
\item For all $t \in \mathbb C,$ except a finite set $T(\xi),$ the intersection $X^{\mathbb C}(\xi, t) \cap \Gamma$  is transversal.
\item The intersection $X^{\mathbb C}(\xi,1) \cap \Gamma_{\infty}$ is transversal.
\end{enumerate}

\end{lemma}

\pf

A few comments on the conditions (1), (2).
Zariski open subsets  sets of $S^{n-1}$ are complements to proper subsets of $S^{n-1}$  given by polynomial equations.
Clearly, any Zariski open subset is a  dense open subset of $S^{n-1}.$
The algebraic variety $\Gamma$ can have singular points, hence transversality conditions (1), (2)  mean that $X^{\mathbb C} (\xi, t)$ transversally intersects any stratum of the Whitney stratification of $\Gamma$ or $\Gamma_{\infty},$ correspondingly.

The proof of Lemma is basically in the line of \cite{Vas1}.  Denote $Reg_{\Gamma}$ the set of affine complex hyperplanes
$X^{\mathbb C} \subset \mathbb C^n$ which are transversal to  $\Gamma.$ It is a Zariski open subset of the complex affine Grassmanian $Gr(n,n-1).$
Correspondingly, the complement  $Tan_{\Gamma}= Gr(n,n-1) \setminus Reg_{\Gamma}$, consisting of all hyperplanes which are {\it not} transversal to $\Gamma$, is Zariski closed.

The Grassmanian $Gr(n,n-1)$ can be parametrized by $X^{\mathbb C}=X(\xi,t), \xi \in \mathbb C^n, \ t \in \mathbb C$ where the parameters
$\lambda \xi, \lambda t, \lambda \in \mathbb C \setminus 0$ correspond the same hyperplane $X^{\mathbb C}.$
Then the Zariski closed set $Tan_{\Gamma}$ can be realized as the intersection
$$Tan_{\Gamma}=\bigcap\limits_{j=1}^J V_j$$
of conical hypersurfaces $V_j$ in $\mathbb C^n \times \mathbb C,$
where $V_j=\alpha_{j}^{-1}(0)$ is the  zero locus of a nonzero homogeneous  polynomial $\alpha_{j}(\xi,t).$

Now introduce the three subsets of the unit sphere $S^{n-1}:$

1.$\Sigma$ - the set of all $\xi \in S^{n-1}$ such that for infinitely many $t \in \mathbb C$ the hyperplane $X^{\mathbb C}(\xi,t)$  {\it does  not} intersect $\Gamma$ transversally.

2. $\Sigma_{\infty}$ - the set of  all $\xi \in S^{n-1}$ such that $X^{\mathbb C}(\xi,1)$ {\it does not} intersect transversally $\Gamma_{\infty}.$

3.  $ W=S^{n-1} \setminus (\Sigma \cup \Sigma_{\infty}).$

Then $W$ possesses properties (1), \ (2) and all we need is  to prove that $W$ contains a Zariski open subset $\Xi \subset S^{n-1}.$

 Fix arbitrary $\xi_0 \in \Sigma$.  By the definition of $\Sigma$, for infinitely many $t,$ the hyperplane $X^{\mathbb C}(\xi, t)$ belongs to the Zariski closed set $Tan_{\Gamma}= Gr(n,n-1) \setminus Reg_{\Gamma}.$ Therefore, as it is pointed out above, there is a finite set of nonzero homogeneous polynomials $\alpha_{j}(\xi,t)$ such that
$$\alpha_{j}(\xi_0,t)=0$$
for infinitely many values of $t.$
Since $t \to \alpha_{j}(\xi_0,t)$ is a polynomial,  $\alpha_{j}(\xi_0, t)=0$ for all $t \in \mathbb C.$
Therefore, if
$$\alpha_{j}(\xi,t)=\sum\limits_{k=0}^N \alpha_{j,k}(\xi)t^k$$
then for each coefficient $\alpha_{j,k}(\xi_0)=0.$
Since $\xi_0 \in S^{n-1}$ is arbitrary, $\Sigma$ is contained in the finite intersection of the zero sets on $S^{n-1}$ of the polynomials $\alpha_{j,k}(\xi):$
$$\Sigma \subset \bigcap\limits_{j,k} \alpha_{j,k}^{-1}(0).$$
The coefficients  $\alpha_{j, k}(\xi)$ are homogeneous polynomials, which are not all identically zero, hence there is $j_0,k_0$ such that the restriction $\alpha_{j_0,k_0}\vert_{S^{n-1}}$ is not zero identically.
Thus, $\Sigma$ is contained in a proper Zariski closed subset $F=\alpha_{j_0,k_0}^{-1}(0)  \subset S^{n-1}.$

The arguments for the set $\Sigma_{\infty}$ (property (2)) are similar. The variety $\Gamma$ is replaced  here by the conical variety (the asymptotic cone) $\Gamma_{\infty}=Q_d^{-1}.$
The set $Reg_{\Gamma_{\infty}} $ of all  $(\xi,t) \in \mathbb C^{n} \times \mathbb C$ such that the hyperplane $X^{\mathbb C} (\xi,t),$ transversally intersects $\Gamma_{\infty},$ is Zariski open. Its complement $Tan_{\Gamma_{\infty}}$
consists of all $(\xi,t)$ such that $X^{\mathbb C}(\xi,t)$ does not transversally intersect the cone $\Gamma_{\infty}$  and is Zariski closed.  Therefore, there exists a nonzero homogeneous polynomial $\beta(\xi,t)$ such that $\beta(\xi,t)=0$  for all $(\xi,t) \in Tan_{\Gamma_{\infty}}.$  Again, due to the homogeneity, $\beta(\xi,t)$ is not identical zero on $S^{n-1} \times \mathbb C.$

Let $\xi \in S^{n-1}$ is not in $ \Sigma_{\infty}.$ By the definition of $\Sigma_{\infty},$ tt means that $(\xi, 1) \in Tan_{\Gamma_{\infty}}. $ Moreover, this is true for all $(\xi,t)$  with $t \neq 0.$ Indeed,
$\Gamma_{\infty}$ is a cone and hence $X^{\mathbb C}(\xi,t), \ t \neq 0$ is not transversal to $\Gamma_{\infty}$ if and only if $(\xi,1)$ does so.  Therefore, if $\xi \notin \Sigma_{\infty}$ then $\beta(\xi, t)=0$ for all $t \neq 0.$ This implies $\beta_j(\xi)=0$ where $\beta(\xi,t)=\sum_{j=0}^B \beta_j(\xi) t^j,$ i.e.
$\Sigma_{\infty}  \subset  \beta_j^{-1}(0)$ for any $j.$ Since $\beta$ is not identical zero on $S^{n-1} \times \mathbb C$ , one of the coefficients $\beta_{j_0}$ is not identical zero on $S^{n-1}$   and $F_{\infty}=\beta_{j_0}^{-1}(0)=F_{\infty}$ is a proper Zariski closed subset of $S^{n-1},$ contained $\Sigma_{\infty}.$

Finally, denote
$$\Xi= S^{n-1} \setminus (F \cup F_{\infty}).$$
Then $\Xi$ is Zariski open subset of $S^{n-1}$ and $\Xi \subset W.$ The latter means that the set $\Xi$ satisfies both properties (1), (2) which completes the proof of Lemma \ref{L:Zariski}.
\qed

\subsection{Projectivization} \label{S:projectivization}

Along with the affine coordinates $z_1, \cdots, z_n$ in $\mathbb C^n,$ we consider the $(n+1)$-dimensional space obtained  by adding the additional coordinate $u_0 \in \mathbb C$ and
introduce the homogeneous coordinates $(u_0; u)= (u_0, u_1, \cdots, u_n)$ by
\begin{equation} \label{E:proj}
z_j= \frac{u_j}{ u_0}, \ j=1, \cdots, n.
\end{equation}
Thus, we pass to the space $\mathbb C^{n+1}$ of the vectors $(u_0, u), u \in \mathbb C^n.$
The space $\mathbb C^n$ can be identified with the hyperplane $u_0=1.$
The hyperplane $u_0=0$ in $\mathbb C^{n+1}_u$ corresponds to the "infinitely distant" hyperplane in the $z$-space.

Sometimes, we will be using the notation
$\mathbb C^{n+1}_u$ to distinguish  the variable $(u_0,u).$  Also, objects in the $u$-space corresponding to related objects in the $z$-space will be often equipped
with the subscript $\widehat{}$ .

Fix $t \in \mathbb C$ and represent it as
\begin{equation} \label{E:t:u}
t=\frac{1}{u_0},
\end{equation}
then the hyperplane $X^{\mathbb C}(\xi,t)=\{ z \in \mathbb C^n: \xi_1 z_1 + \cdots \xi_n z_n=t \}$  transforms into the $(n-1)$-plane in the space $\mathbb C^{n+1}_u:$
\begin{equation} \label{E:widehatX}
\widehat X^{\mathbb C}(\xi,u_0) = \{u_0\} \times X^{\mathbb C}(\xi,1)=  \{ (u_0, u): u \in \mathbb C^{n} : \xi_1 u_1+ \cdots +\xi_n u_n=1 \}.
\end{equation}

We also denote $\widehat X^{\mathbb C}(\xi)$ the $n$-plane in $\mathbb C^{n+1}_u:$
\begin{equation} \label{E:widehatXksi}
\widehat X^{\mathbb C}(\xi)=\mathbb C \times X^{\mathbb C}(\xi,1)=
\{ (u_0, u) \in  \mathbb C \times \mathbb C^n=\mathbb C^{n+1}: \xi_1 u_1+...\xi_n u_n=1 \},
\end{equation}
which can be foliated into $(n-1)$-planes as
\begin{equation} \label{E:union}
\widehat X^{\mathbb C} (\xi)= \bigcup\limits_{u_0 \in \mathbb C} \widehat X^{\mathbb C}(\xi,u_0).
\end{equation}
Let
$$Q=Q_d+Q_{d-1}+ \cdots + Q_0, \ d=deg Q,$$
be the  decomposition of the polynomial $Q$ into the sum of homogeneous polynomials.
By the definition
$$\Gamma=\{z \in \mathbb C^n: Q(z)=Q_d(z)+ \cdots + Q_0(z)=0 \} .$$

Fix $u_0 \neq 0.$ Substitution $z=\frac{u}{u_0}$ transforms the equation $Q(z)=0$ into equation
\begin{equation} \label{E:Q}
\widehat Q(u_0,u)= Q_d(u)+u_0 Q_{d-1}(u) + \cdots + u_0^d Q_0(z)=0.
\end{equation}
Finally, let us introduce the $n$-dimensional algebraic variety  $\widehat \Gamma \subset \mathbb C^{n+1}:$
\begin{equation} \label{E:widehatG}
\widehat \Gamma=\{(u_0,u) \in \mathbb C \times \mathbb C^n  : \widehat Q(u_0,u)=0 \}.
\end{equation}

For fixed $u_0 \in \mathbb C$ denote
$\widehat \Gamma_{u_0}=\{ u \in \mathbb C^n: Q(u_0, u)=0 \}.$
The value $u_0=1$ corresponds to the original variety $\Gamma=Q^{-1}(0)$ (the complexification of $\partial K$):
\begin{equation} \label{E:G1}
\widehat \Gamma_1= \{u_0=1\} \cap \widehat \Gamma = \{1\} \times \Gamma.
\end{equation}
 The exceptional value $u_0=0$ defines the cross-section of $\widehat \Gamma$ by the "infinitely distant" plane $u_0=0$ (which corresponds to $t=\frac{1}{u_0}=\infty$) and is realized as the asymptotic cone  $\Gamma_{\infty}$ which is the zero set of the leading homogeneous term $Q_d$ in the decomposition of the polynomial $Q:$
\begin{equation}\label{E:G0}
\widehat \Gamma_0 = \{0\} \times \Gamma_{\infty}= \{0\} \times Q_d^{-1}(0).
\end{equation}

\section{Fiber bundle}

Fix a unit vector $\xi \in S^{n-1}$ and consider  the intersection of the $n$-dimensional plane $\widehat X^{\mathbb C}(\xi)$ in $\mathbb C^{n+1}_u$  (\ref{E:widehatXksi}) and the $n$-dimensional projective variety  $\widehat \Gamma$ (\ref{E:widehatG}). The foliation (\ref{E:union}) of the hyperplane
$\widehat X^{\mathbb C}(\xi)$ induces the foliation of its intersection  with $\widehat \Gamma$ :
\begin{equation} \label{E:foli}
\widehat X^{\mathbb C}(\xi) \cap \widehat \Gamma = \bigcup\limits_{u_0 \in \mathbb C} (\widehat X^{\mathbb C}(\xi, u_0) \cap \widehat \Gamma).
\end{equation}

We will be using the following projective version of Lemma \ref{L:Zariski}:
\begin{lemma} \label{L:Zariskihat}
Let $\Xi$ be the Zariski open set from Lemma \ref{L:Zariski} and let $\xi \in \Xi.$ Then
there is a finite set $\widehat{T} (\xi) \subset \mathbb C$ such that
\begin{enumerate}
\item For all $u_0 \notin \widehat{T} (\xi)$ the hyperplane $\widehat X(\xi, u_0)$ intersects $\widehat \Gamma$ transversally.
\item  $0 \notin \widehat {T} (\xi).$
\end{enumerate}
\end{lemma}

\pf
By the construction, the exceptional sets in Lemmas \ref{L:Zariski} and \ref{L:Zariskihat} are related to each other by
$$T(\xi)= \{ t =\frac{1}{u_0}: u_0 \in \widehat T(\xi), u_0 \neq 0 \},$$
$u_0=0$ corresponds to $t = \infty.$ Property (1) in Lemma \ref{L:Zariski} transforms into property (1) in Lemma \ref{L:Zariskihat}. Property (2) in Lemma \ref{L:Zariski}  claims that $X(\xi,1)$ transversally intersects the asymptotic cone $\Gamma_{\infty}$ which implies that the infinitely distant hyperplane $u_0=0$ intersects transversally $\widehat \Gamma.$ This means that $u_0=0$ is not in the exceptional set, i.e., $0 \notin \widehat {T}(\xi)$ which is just property (2) in Lemma \ref{L:Zariskihat}.
\qed

Lemma \ref{L:Zariskihat}  implies that if $\xi \in \Xi$ then the foliation (\ref{E:foli}) is in fact a regular fibration out of the finite singular set $\widehat T(\xi)$ where the intersections with the planes are not transversal.
Now we want to restrict this fibration on a certain domain in the punctured complex plane $\mathbb C \setminus \widehat T(\xi).$

\begin{lemma} \label{L:fiber}
Consider the projection mapping
$$\pi_{\xi}: \widehat X^{\mathbb C} (\xi) \cap \widehat \Gamma \to \mathbb C,$$
defined by:
\begin{equation}\label{E:pi}
\pi_{\xi}(u_0, u)=u_0.
\end{equation}
Let $\Omega $ be a domain in the complex plane satisfying the two conditions:
\begin{enumerate}
\item $\Omega \subset \mathbb C \setminus \widehat T(\xi)$,
\item $\mathbb C \setminus \Omega$ is a connected set.
\end{enumerate}
Consider the restriction
\begin{equation} \label{E:restr}
\pi_{\xi, \Omega} : \pi_{\xi}^{-1}(\Omega ) \to \Omega
\end{equation}
of the projection $\pi_{\xi}$ onto the pre-image of the domain $\Omega.$
Then the  mapping $\pi_{\xi,\Omega}$ defines a trivial fiber bundle over the base $\Omega.$
\end{lemma}

\pf  The $\pi_{\xi}$- fibers over $u_0 \in \mathbb C$ are
$$\pi_{\xi}^{-1}(u_0)= \widehat X^{\mathbb C}(\xi, u_0) \cap \widehat \Gamma.$$
When $u_0 \notin \widehat T(\xi),$ then the intersections are transversal (Lemma \ref{L:Zariskihat}). Therefore, by Tom's first isotopy lemma, see e.g. \cite{T, G}, (\cite{JM}, Prop. 11.1).
the mapping $\pi_{\xi}$ is a locally trivial fiber bundle over the punctured complex plane $B=\mathbb C \setminus \widehat T(\xi).$

The restricted fiber bundle $\pi_{\xi,\Omega}$ over the base $\Omega$ is also locally trivial since $\Omega \subset B$ due to 1). The base $\Omega$ of this fiber bundle is contractible due to condition 2). The locally trivial fiber bundle over a contractible base
 is homotopically fiberwise equivalent to a trivial fibration (\cite{Sp},
 Ch. 2, S. 8, Cor. 15)., i.e. to the direct product of $\Omega$ and any fixed  single fiber $\pi_{\xi}^{-1}(\zeta).$
\qed

%\section{Analytic extension of $A_K(\xi,t)$}

%\subsection{Integration on chains}

\section{Integration on chains}

We introduce the orientation on the hyperplane $X(\xi,t) \subset \mathbb R^n$ induced by the normal vector $\xi.$
Namely, the basic $\tau_1, \cdots, \tau_{n-1}$ in $X(\xi,t)$ is right-oriented if \\
$det (\xi, \tau_1, \cdots, \tau_{n-1}) >0 . $ Then the cross-sections $c=X(\xi,t) \cap K$  can be considered as $(n-1)$-dimensional singular
chains with the boundaries $\partial c \subset \partial K.$
The value $A_K(\xi,t)$ of the $(n-1)$-dimensional volume of the intersection $K \cap X(\xi,t)$ can be expressed as the integral on such chains:
\begin{equation} \label{E:AKu}
A_K(\xi,t)=\int\limits_{c} \omega_{\xi}(x), \  c= X(\xi,t) \cap K,
\end{equation}
against
the differential $(n-1)$-form
\begin{equation}\label{E:omega}
\omega_{\xi}(x)= \sum\limits_{j=1}^n (-1)^{j+1}\xi_j dx[j],
\end{equation}
$dx[j]=dx_1 \wedge \cdots \wedge dx_{j-1} \wedge dx_{j+1} \wedge \cdots \wedge dx_n.$

The differential form $\omega_{\xi}(x)$ extends to $\mathbb C^n$ as a  holomorphic $(n-1)$- differential form  $\omega_{\xi}(z)$ by writing $dz_j$ in place of $dx_j.$  the Then the integral $\int_c \omega_{\xi}(z)$ is well defined for
any singular chain $c \subset \mathbb C^{n}$ of the real dimension $n-1.$

Let us transfer now the above construction  to the projective model described in Section \ref{S:projectivization}. To underline the passage from the $x$- or $z$-space to the $u$-space, we will be using for the chains of integration in the projective model the subscript $\widehat {} $ \ .

Firstly, let us look at the real space. For real $u_0,$ the mapping
$$ (u_0, u) \to z=\frac{u}{u_0}$$
Mapping (\ref{E:proj}) transforms the domain $K \subset \mathbb R^n$ into the domain
$$ \widehat {u_0 K} = \{ u_0, u_0 x): x \in K \} =\{u_0\} \times u_0 K.$$
The real $(n-1)$-plane $X(\xi,t)$ in $\mathbb R^n$  becomes the real $(n-1)$-plane $\widehat X(\xi,u_0)= \{ (u_0,u): \xi_1 u_1 + \cdots + \xi_n u_n =1 \}, \ t =\frac{1}{u_0},$ in $\mathbb R^{n+1}$ (see (\ref{E:widehatX})).

Therefore, the volume integral (\ref{E:AKu}) transforms as:
\begin{equation}\label{E:AKupro}
A_K(\xi,t)=\frac{1}{u_0^{n-1}} \int\limits_{ \widehat c(\xi,u_0)} \omega_{\xi}(u),
\end{equation}
where $u_0$ and $t$ are real, $t=\frac{1}{u_0},$ and the integral is taken over the $(n-1)$-dimensional
chain $\widehat c(\xi, u_0)$ in the real space $R^{n+1}:$
\begin{equation} \label{E:chat}
\widehat c(\xi,u_0)= \{ (u_0, u) \in \mathbb R^{n+1}: u \in u_0 K, \sum\limits_{j=1}^n \xi_j u_j =1\}.
\end{equation}
The orientation on $\widehat c(\xi, u_0)$ is induced by the orientation of the hyperplane $X(\xi,1) \subset \mathbb R^n.$

The integral in the right hand side of (\ref{E:AKupro}) extends as a function $I_{|xi}$ defined on
 singular chains $\widehat c$ of real dimension $n-1,$  belonging to the complex space  $\mathbb C^{n+1},$ by:
\begin{equation} \label{E:Ic}
I_{\xi}(\widehat c)=\int\limits_{\widehat c} \omega_{\xi}(u).
\end{equation}

By (\ref{E:AKu}) the real values  $A_K(\xi,t) >0 $ of the section function can be computed in terms of the above integral as
\begin{equation}\label{E:AKupro}
A_K(\xi,t)= u_0^{n-1} I_{\xi}( \widehat c(\xi, u_0)), \ t=\frac{1}{u_0},
\end{equation}
where $u_0$ is real, $u_0 \neq 0,$ the chain of integration $\widehat c(\xi,u_0)$ is defined by (\ref{E:chat}). The boundary of the chain $\widehat c(\xi, u_0)$ belongs to the real hypersurface $\partial (u_0 K) \cap (\{u_0\} \times X(\xi, t)), \ t=\frac{1}{u_0}, $ which, in turn,  is contained in the real part of
$\widehat X^{\mathbb C}(\xi, u_0)  \cap \widehat \Gamma.$

\begin{lemma} \label{L:c1c2} Let $\Xi \subset S^{n-1}$ be the Zariski open set from Lemma \ref{L:Zariskihat}. Let $\xi \in \Xi$  and $u_0 \notin \widehat T(\xi),$
where $\widehat T(\xi)$ is the exceptional set in Lemma \ref{L:Zariskihat}.
 Let $\widehat c_1, \widehat c_2$ be two singular chains in $\mathbb C^{n+1}$ of the real dimension $n-1$ which are homologically equivalent (
 $\widehat c_1 \sim \widehat c_2$)  relatively
$ \widehat X^{\mathbb C}(\xi, u_0)  \cap  \widehat \Gamma,$ which means that
\begin{enumerate}
\item $\partial \widehat c_1, \partial \widehat c_2 \subset \widehat X^{\mathbb C}(\xi, u_0)  \cap  \widehat \Gamma,$
\item there exists a $(n-1)$-chain $\widehat d \subset  X^{\mathbb C}(\xi, u_0) \cap \widehat \Gamma$ and a $n$-chain $\widehat e$ such that
$\widehat c_1- \widehat c_2 + \widehat d=\partial \widehat  e .$
\end{enumerate}
Then $I_{\xi}(\widehat c_1)=I_{\xi}(\widehat c_2).$
\end{lemma}

\pf
Since $d\omega_{\xi}=0,$ Stokes formula yields
$$I_{\xi}(\widehat c_1)-I_{\xi}(\widehat c_2)+I_{\xi} (\widehat d)=\int_{\partial \widehat e} \omega_{\xi}(u)=\int_{\widehat e}d\omega_{\xi}(u)=0.$$
Due to the choice of $\xi$ and $u_0,$ the intersection $X^{\mathbb C}(\xi, u_0) \cap \widehat \Gamma,$ is transversal. This is an algebraic variety of the complex dimension $n-2.$  Since the chain $e$ is contained in this variety, the holomorphic $(n-1)-$ differential form $\omega_{\xi}(u)$ integrates over $\widehat d$ to $0.$
Therefore, $I_{\xi}(\widehat c_1)-I_{\xi}(\widehat c_2)=0.$
\qed

Denote
$$H_{n-1}(\xi,u_0)= H_{n-1}(\mathbb C^{n+1}, \widehat X^{\mathbb C}(\xi,u_0) \cap \widehat \Gamma)$$
the $(n-1)$-dimensional homology group, relative with respect to $  \widehat X^{\mathbb C}(\xi, u_0) \cap \widehat \Gamma.$
The elements of this group are the equivalence classes in the space of the $(n-1)$-chains $\widehat c$ satisfying condition (1) in Lemma \ref{L:c1c2},
with respect to the equivalence relation given by condition (2) in the same lemma.

Lemma \ref{L:c1c2} states that the function $I_{\xi}(\widehat c)$ takes same values on the relatively equivalent chains $\widehat c_1 \sim \widehat c_2$
and hence it is defined on the elements of the relative homology groups $H_{n-1}(\xi,u_0)$. This proves the following
\begin{lemma}\label{L:L}
Let $\xi \in \Xi.$ The equality
\begin{equation} \label{E:L}
L_{\xi}([\widehat c])=I_{\xi}(\widehat c),
\end{equation}
defines, for $ u_0 \notin \widehat T(\xi),$ a linear function
$$L_{\xi}: H_{n-1}(\xi, u_0) \to \mathbb C,$$
where $[\widehat c] \in H_{n-1}(\xi, u_0)$ the relative homology class containing the $(n-1)$-chain $\widehat c.$
\end{lemma}

\section{Singe-valued analytic extension of the section function and its growth at $\infty$}

\subsection{The strategy}
Equalities (\ref{E:AKu}), \ (\ref{E:AKupro}) show that the value $A_K(\xi,t), \ t \neq 0,$ can be computed (as long as $X(\xi,t)$ meets $K$ and is transversal to $\partial K$, i.e., when $A_K(\xi,t)>0$)  via the integral
\begin{equation} \label{E:multi}
A_K(\xi,t)= \frac{1}{u_0^{n-1}}  L_{\xi} ([\widehat c_{\xi, u_0}]),
\end{equation}
where $ u_0=\frac {1}{t}$ and the chain $\widehat c(\xi, u_0)$ in the $u$-space is generated by the chain $X(\xi,t) \cap K$ in the $x$-space (formula (\ref{E:chat})).

While, for given $\xi \in S^{n-1}$ and for real $t,$ such that $A_K(\xi,t) >0 ,$ the choice  of the chain $c(\xi,t)$ in formula (\ref{E:multi}) is determined, for complex $t \notin T(\xi),$  or, equivalently, for complex $\frac{1}{t}=u_0 \notin \widehat T(\xi),$ the linear function $L_{\xi}$ is defined on homological classes (Lemma \ref{L:L}) and hence can be viewed as a multi-valued function of $t.$  Therefore, to construct an analytic continuation of $A_K(\xi,t)$ into punctured complex plane we need to construct a single-valued analytic branch, i.e., a continuous section of the corresponding homological vector bundle. The initial value for this section
can be chosen by fixing a real cross-section $X(\xi, t_0) \cap \partial K$ and demanding that the analytic extension coincides with $A_K(\xi, t)$ for real $t$ in a neighborhood of $t_0.$

The above construction is convenient to realize in the projective model, i. e,, in the $u$-space.

\subsection{Homological vector bundle}

Fix a vector $\xi$ in the Zariski open set $\Xi \subset S^{n-1}$ from Lemma \ref{L:Zariski}. By the construction, the set of all $t \in \mathbb C$ such that the hyperplane $X(\xi,t)$ intersects
$ \Gamma$ non-transversally is finite. Chose real $t_0$ which possesses such a property and $X(\xi,t_0)$ meets $\partial K$ transversally. Without loss of generality we can assume that $t_0 =1.$ Then $u_0=\frac{1}{t_0}=1$ does not belong to the set $\widehat T(\xi)=\{\frac{1}{t}: t \in T(\xi) \}$ of exceptional values of $u_0.$ Also, by Lemma \ref{L:Zariskihat}, $0 \notin \widehat T(\xi).$

Let $\Omega \subset \mathbb C \setminus \widehat T(\xi)$  be a complex domain, as in Lemma \ref{L:fiber}. Since $0, 1 \notin \widehat T(\xi),$ we can choose
$\Omega$ so that $0, 1 \in \Omega.$

Consider the fiber bundle with the base space $\Omega$, constructed in Lemma \ref{L:fiber}:
$$\pi_{\xi}: \widehat X^{\mathbb C} (\xi) \cap \widehat \Gamma \to \mathbb C,$$
and its restriction
$$\pi_{\xi,\Omega}: \pi_{\xi}^{-1}(\Omega) \to \Omega$$
on the set $\Omega.$
Lemma \ref{L:fiber} claims that the mapping $\pi_{\xi,\Omega}$ defines a trivial fiber bundle with the base $\Omega$ and the fibers $\widehat X^{\mathbb C}(\xi, u_0),  \ u_0 \in \Omega.$

Consider the associated vector bundle
of the relative homology groups $H_{n-1}(\xi, u_0)=  H_{n-1}(\mathbb C^{n+1}, \widehat X^{\mathbb C}(\xi, u_0) \cap \widehat \Gamma).$
The fiber over the point $u_0 \in \Omega$ is the pair $(u_0, H_{n-1}(\xi,u_0)).$

Denote
$$\widetilde \pi_{\xi} : \bigcup\limits_{u_0 \in \mathbb C} (u_0, H_{n-1}(\xi, u_0))  \to \mathbb C$$
the projection mapping: $\widetilde \pi (u_0, H_{n-1}(\xi, u_0))=u_0$
and
$$\widetilde \pi_{\xi, \Omega} : \widetilde \pi_{\xi}^{-1}(\Omega) \to \Omega$$
the restriction $\widetilde \pi_{\xi}$ to the pre-image of the  $\Omega.$

Since the restricted fiber bundle $\pi_{\xi,\Omega}$ over $\Omega$ is trivial (Lemma \ref{L:fiber}) , the associated homological  vector bundle $\widetilde \pi_{\xi,\Omega}$  is also trivial.
The trivial vector bundle admits a continuous section (the right inverse mapping) $\widetilde \pi_{\xi, \Omega} ^{-1}$
$$\widetilde \pi_{\xi,\Omega}^{-1} : \Omega \to H_{n-1}(\xi,\Omega)=  \bigcup \limits_{u_0 \in \Omega} H_{n-1}(\xi,u_0),$$
with the prescribed value
$$\widetilde \pi_{\xi,\Omega}(1)=[\widehat c_{\xi,1}]$$
where $\widehat c_{\xi,1}$ is the $(n-1)$-chain (\ref{E:chat}) produced by the cross-section $X(\xi,1) \cap K$ in the original real $x$-space.

\subsection{The analytic extension of $A_K(\xi,t)$ for complex $t$}

Since now, the coordinate $u_0$ plays a role of a variable in the complex domain $\Omega,$  hence it would be natural to change the notation
$$u_0 =\zeta.$$
The value of the mapping $\widetilde \pi_{\xi,\Omega}(\zeta)$  at the point $\zeta \in \Omega$ is an element ( relative homology class) in $H_{n-1}(\xi,\zeta),$  continuously varying when $\zeta$ runs along the path $C \subset \Omega$ from $1$ to $0.$

Finally, consider the  composition $F$ of the "lifting" mapping $\widetilde \pi^{-1}_{\xi,\Omega}$ of $\Omega$ to the homology groups
and the integral mapping $L_{\xi}$  (\ref{E:L}) of the homologies to the complex plane:
$$F_{\xi} (\zeta)= L_{\xi} (  \widetilde \pi^{-1}_{\xi,\Omega}(\zeta)) .$$  The following diagram is commutative:

 \begin{center}
\begin{tikzpicture}[every node/.style={midway}]
\matrix[column sep={4em,between origins},
        row sep={2em}] at (0,0)
{ \node(H)   {$H_{n-1,\Omega}$}  ; & \node(C) {$\mathbb C $}; \\
  \node(Omega) {$\Omega$};                   \\};
\draw[->] (H) -- (Omega) node[anchor=east]  {$\widetilde\pi_{\xi, \Omega}$};
%\draw[<-] (H) -- (Omega) node[anchor=east]  {$\widetilde\pi^{-1}_{\xi, \Omega}$};
\draw[->] (Omega) -- (C) node[anchor=north]  {$F_{\xi}$};
\draw[->] (H)   -- (C) node[anchor=south] {$L_{\xi}$};
\end{tikzpicture}
\end{center}

By the construction, $F_{\xi}(\zeta)$ is a single-valued analytic function in the complex domain $\Omega.$ The value $\zeta=u_0=1$
 corresponds to the cross-section  of $\partial K$  by the real hyperplane $X(\xi,1),$  the value $\zeta=u_0=0$ corresponds in the original affine representation to the intersection  $\Gamma$ by the complex hyperplane $X^{\mathbb C}(\xi,t)$ at $t=\infty $ (\ref{E:G1}, \ref{E:G0}).

The variables $\zeta= u_0$ and $t$ are related by $t=\frac{1}{\zeta}.$  Domain $\Omega$ contains the points $0$ and $1,$ and
when $t$ is real and belongs to a neighborhood of $t_0=1$ then $\zeta$ is also real and in a neighborhood of $1,$ and $A_K(\xi,t)$ and $F_{\xi}(\zeta)$ are related, due to (\ref{E:multi}), by
\begin{equation}\label{E:FA}
A_K(\xi, t)= \frac{1}{\zeta^{n-1}} F_{\xi}(\zeta)= t^{n-1} F_{\xi}(\frac{1}{t}).
\end{equation}

\subsection{The upper bound for $\deg p(\xi,t)$}

\begin{lemma} \label{L:degree} For any $\xi \in S^{n-1},$  the degree of the polynomial $ t \to p(\xi,t)$ is at most $m(n-1).$
\end{lemma}
\pf
Fix $\xi \in \Xi.$
The relation (\ref{E:FA}),  and the basic equation (\ref{E:equation}) ) imply
$$p(\xi,\frac{1}{\zeta} )= - q(\xi) A_K^m(\xi, \frac{1}{\zeta}) = - q(\xi) \frac{1}{\zeta^{m(n-1)}} F_{\xi}^m(\xi, \zeta),$$
when $\zeta$ is a real number in a neighborhood of $\zeta_0=1.$

By the uniqueness theorem , the equality extends to all $ \zeta$ in the domain  $\Omega,$ where $F_{\xi}$ is analytic.
Since $0 \in \Omega,$ the function $F_{\xi}$ takes a finite value $F_{\xi}(0)$ at $\zeta=0$ and we have
$$ \zeta^{n-1} \lim\limits_{\zeta \to 0, \zeta \in \Omega} p(\xi, \frac{1}{\zeta})= - q(\xi) F^m_{\xi}(0). $$

The open set $\Omega$ in Lemma \ref{L:fiber} is connected and contains the points $0$ and $1,$ therefore $\Omega$
contains a continuous path $C$ joining $0$ and $1.$  The above limit relation implies, after substitution $\zeta=\frac{1}{t},$ the finiteness of the limit
$$ \lim\limits_{ t \to \infty, \ t \in \widetilde C} \frac{p(\xi,t)}{t^{m(n-1)}} = - q(\xi)F^m_{\xi}(0), $$
where the path $$ \widetilde C=\{t=\frac{1}{\zeta}: \zeta \in C \} $$
joins $1$ and $\infty.$
Thus, the limit in the right hand side is finite and therefore the degree of the polynomial $t \to p(\xi,t)$ does not exceed $m(n-1).$ This is true for $\xi$ in the Zariski open set $\Xi.$
Since $\Xi$ is dense in $S^{n-1},$ it follows, by continuity with respect to $\xi,$ that $\deg_t p(\xi,t) \leq m(n-1)$ for any  vector $\xi \in S^{n-1}.$
\qed

\section{Zeros of  polynomial $p(\xi,t)$}

Denote
$$h_k(\xi)=\sup\limits_{ x \in K} x \cdot \xi$$
the support function of the body $K.$
Given a point $a \in \partial K,$ the equation of the tangent plane to $\partial K$ at $
a$ is
$$T_a(\partial K)=\{ x \in \mathbb R^n: x \cdot \nu_a =h_K(\nu_a),$$
where $\nu_a$ is the external unit normal vector to $\partial K.$ In our previous notations, $T_a(\partial K)=X(\xi,t)$ with $\xi=\nu_a$ and $t=h_K(\xi).$

The point $a \in \partial K$ is {\it elliptic} if the second fundamental form at this point is positive definite. After a suitable rotation, we can assume
that $\nu_a=(0, \cdots, 0, 1)$ and then the convex hypersurface $\partial K$ near the point $a$ can be represented as the graph
$x_n=\psi(x^{\prime}), \  x^{\prime}=(x_1, \cdots, x_{n-1}),$ of a $C^{\infty}$ function

\begin{equation}\label{E:psi}
\psi(x_1,\cdots,x_{n-1})=a_n - \frac{1}{2} \sum\limits_{j=1}^{n-1} \kappa_j (x_j-a_j)^2 + o(|x^{\prime}-a^{\prime}|), \ x^{\prime} \to a^{\prime}.
\end{equation}
Here $\kappa_j $ are the principal curvatures at $a$ and all $\kappa_j >0 $ at elliptic points.

\begin{lemma} \label{L:a}
Let $a \in \partial K$ be an elliptic point and $\xi=\nu_a$
the external normal vector.
Then
\begin{equation} \label{E:t-t0}
A_K(\xi,t)= c (t_0-t)^{\frac{n-1}{2}} (1+ o(t_0-t)), \ t \to t_0-,
\end{equation}
where $t_0=h_K (\xi).$
\end{lemma}

\pf By applying rotation and translation, it suffices to prove Lemma for the case $a=0$ and $\xi=(0,\cdots,0, 1)$ and $h_K(\xi)=t_0=0.$
By (\ref{E:psi}), then $\partial K$ near $a$ writes as
$$ x_n = - \frac{1}{2}(\sum\limits_{j=1}^{n-1} \kappa_j x_j^2 )(1+o (|x^{\prime}| )), \  x^{\prime} \to 0.$$
Then the volume of the section the body $K$ by the hyperplane $x_n=t, \ t>0$ coincides , up to $o(t),$ with the volume of an ellipsoid:
$$A_K(\xi,t)= c_n \frac{1}{\kappa_1 \cdots \kappa_{n-1}} (-t)^{\frac{n-1}{2}}( 1+o(t)), \ t \to 0-,$$  where $c_n$ is the area of the unit sphere $S^{n-1}.$
This proves Lemma.
\qed

\begin{lemma}
The polynomial $p(\xi,t)$ has the following representation
\begin{equation}\label{E:represents}
p(\xi,t)=C(\xi) (h_K(\xi)-t)^{\frac{m(n-1)} {2} } (h_K(-\xi)+t)^{\frac{m(n-1)}{2}} .
\end{equation}

\end{lemma}

\pf
The normal vectors $\xi=\nu_a$ at elliptic points $a \in \partial K$ constitute a dense subset of the unit sphere $S^{n-1}.$ Indeed, since $\partial K$ is convex, the second fundamental form is non-negative. Therefore, non-elliptic points are those points on $\partial K$ at which at least one principal curvature $\kappa_j=0,$ i.e.,
Gaussian curvature $\kappa(a) =\kappa_1 \cdots \kappa_{n-1}=0.$ Gaussian curvature $\kappa (a)$ is the determinant of the Gauss mapping
 $\gamma: \partial K \to S^{n-1}, \ \gamma(a)=\nu_a, a \in \partial K$ and therefore the vectors $\xi=\nu_a$ at the non-elliptic points $a \in \partial K$ constitute the
set of critical values of the mapping $\gamma.$ By Sard theorem this set has the measure zero.

Thus, the normal vectors $\xi=\nu_a$ at elliptic points are dense on $S^{n-1}.$ By continuity, it follows that for all $\xi \in S^{n-1}$ the polynomial $p(\xi,t)$ vanishes at $t=h_k(\xi)$ to the order $\frac{m(n-1)}{2}.$ In particular,relation (\ref{E:t-t0}) holds for all $\xi \in S^{n-1}.$

Fix $\xi \in S^{n-1}.$ The function $t \to A_K(\xi, t )$  has zero at $t_+=h_K(\xi)$ of the order $\frac{m(n-1)} {2}.$ Correspondingly, the function $A_K(-\xi, t)$
has zero at $t_=h_K(\-xi)$ of the same order.
From the evenness property of Radon transform: $A_K(\xi,-t)=A_K(-\xi, t),$ we conclude that $A_K(\xi, t)$ has two zeros: $t_+=h_K(\xi)$ and $t_=-h_K(-\xi),$ each one of order at least $\frac{m(n-1)}{2}).$  The equation $q(\xi)A_K^m(\xi,t)= -p (\xi,t)$ implies, due to Bezout theorem, that the polynomial $p(\xi,\cdot)$ is divisible both
by $(t_{+}-t)^{\frac{m(n-1)}{2}}$ and by $(t-t_{-} )^{\frac{m(n-1)}{2} } .$ Since the total degree does not exceed $m(n-1)$ (Lemma \ref{L:degree}), the polynomial $p(\xi, \cdot)$ has the  form (\ref{E:represents}).
\qed

\section{End of the proof of Theorem \ref {T:main}}

Since $$A_K(\xi,t)=\sqrt [m] { \frac{-p(\xi,t)}{q(\xi)} },$$ representation (\ref{E:represents}) implies
$$A_K(\xi,t)=C(\xi) (h_K(\xi)-t)^{\frac{n-1}{2}} (h_K(-\xi)+t)^{\frac{n-1}{2}},$$
for $-h_K(-\xi) <t < h_K(\xi),$ where $C(\xi)$ is a new coefficient depending on $\xi.$
The interval $t \in (-h_K( -\xi), h_K(\xi)) $   corresponds to the hyperplanes  $X(\xi,t)$ which meet the interior of $K,$ i.e.,  $A_K(\xi,t) >0.$

\subsection{Moment conditions}
The section function is the Radon transform $A_K=\mathcal R(\chi_K)$ of the characteristic function of the body $K,$  therefore it satisfies the range condition for Radon transform \cite{He}.
Namely, the $k$-th power moment
$$M_k(\xi)=\int_{-\infty}^{\infty}A_K(\xi,t) t^k dt=
C(\xi) \int\limits_{-h_K(-\xi)}^{h_K(\xi)} (h_K(\xi)-t)^{\frac{m(n-1)}{2}} (t+h_K(-\xi))^{\frac{m(n-1)}{2}}  t^k dt$$
extends from the sphere $|\xi|=1$ to $\mathbb R^n$ as homogeneous polynomial of degree $k.$

Denote for simplicity, $a=a(\xi)=-h_K(-\xi), \ b=b(\xi)=h_K(\xi)$ and change the variable in the integral
$$s=\frac{2t- a- b }{b-a}.$$
Then
$$M_k(\xi)=C(\xi) \int\limits_{-1}^{1} (1-s^2)^{\frac{n-1}{2}}( s(b-a)+(a+b))^k ds,$$
where $C(\xi)$ is a new function of $\xi.$

We are interested in the first three moments $M_0(\xi), M_1(\xi), M_2(\xi).$
We have
$$M_0(\xi)= C(\xi) \int\limits_{-1}^{1} (1-s^2)^{\frac{n-1}{2}} ds,$$
By the moment condition, $M_0(\xi)$ is (a nonzero) constant, $M_0(\xi)=M_0$ and hence $C(\xi)=C_0$ is constant.

The first power moment is
$$M_1(\xi)=C_0 (a(\xi)+b(\xi)) M_0, $$
because $\int\limits_{-1}^{1} (1-s^2)^{\frac{n-1}{2}} s ds=0.$ Therefore, $a(\xi)+b(\xi)=h_K(\xi)-h_K(-\xi)$ coincides on  $|\xi|=1$ with a homogeneous linear polynomial
\begin{equation}\label{E:hk}
h_K(\xi)- h_K(-\xi)= e \cdot \xi ,
\end{equation}
where $e$ is a vector in $ \mathbb R^n.$

Now,
$$M_2(\xi)=C_0 (C_1(b-a)^2 + M_0 (a+b)^2) .$$
where we  have denoted $C_1=\int\limits_{-1}^{1} (1-s^2)^{\frac{n-1}{2}} s^2 ds.$

By the moment conditions,  $M_2(\xi)$ extends form the unit sphere as a homogeneous quadratic polynomial.  Since $(a(\xi)+b(\xi))^2= ( e \cdot \xi )^2$ also extends as a
quadratic polynomial, we conclude that $(b(\xi)-a(\xi))^2=(h_K(\xi)+h_K(-\xi))^2$   does so.

\subsection{Proving that $\partial K$ is an ellipsoid}
Denote
$$H(\xi)=h_K(\xi)- \frac{1}{2} e \cdot \xi.$$
Then we have from (\ref{E:hk}):
$$H(\xi)-H(-\xi)=0, \  H(\xi)=\frac{1}{2} (H(\xi)+H(-\xi)) = \frac{1}{2} ( h_K(\xi) + h_K(-\xi)).$$
and hence $ H^2(\xi)= \frac{1}{4} (h_K(\xi)+h_K(-\xi))^2$ extends as a quadratic polynomial. But $H(\xi)$ is the support function of the translation $K_1=K - \frac{1}{2} e :$
$$H(\xi)=h_K(\xi) - \frac {1}{2} e \cdot \xi= h_{K_1}(\xi).$$
 Thus, the support function $h_{K_1}$ is the square root of a quadratic polynomial and this implies that $\partial K_1$ is an ellipsoid. Then the translated set $\partial K=\partial K_1 -\frac{1}{2}e$ is also an ellipsoid.
Theorem \ref{T:main} is proved.

\section{Concluding remarks}
\begin{itemize}
\item It was proved in \cite{A1} that the polynomial integrability of a body with infinitely smooth boundary implies the convexity of the body.
Possibly, also the convexity condition in Theorem \ref{T:main} is redundant.
\item The question considered in this article  can be also addressed unbounded domains $K.$  The set of second order convex hypersurfaces in $\mathbb R^n$
 consists of the family of ellipsoids and of two families of unbounded hypersurfaces: a single sheet of two-sheets elliptic hyperboloids and elliptic paraboloids.
 The computation shows that the section function of a convex domain, whose boundary is one of the above listed quadratic surfaces, satisfies equation (\ref{E:equation}) for $0 <A_K(\xi,t) < \infty,$ with $m=2.$  It might be of interest to explore whether  there are other solutions of (\ref{E:equation})? Theorem \ref{T:main} answers this question (in negative) for the case of bounded domains.
\end{itemize}

%    Insert the bibliography data here.

\noindent
Bar-Ilan University and Holon Institute of Technology;  Israel.

\noindent
{\it E-mail address:} agranovs@math.biu.ac.il


\begin{thebibliography}{15}
\bibitem {A1} M. Agranovsky, {\em On polynomially integrable domains in Euclidean spaces,} in Complex Analysis and Dynamical Systems, Trends in Mathematics, Birkhauser, Springer, 2018; arXiv.1701.05551.
\bibitem{A2} M. Agranovsky, {\em On algebraically integrable domains}, Contemp. Math., vol.733, 33-44;   arXiv.1705.06063v2.
\bibitem{A3} M. Agranovsky, { Locally polynomially integrable surfaces and finite stationary phase expansions}, J. d'Analyse Math., 141 (2020), pp. 23-47.
\bibitem{A4} M. Agranovsky, {\em Domains with algebraic $X-$ray transform}, Anal. Math. Phys., 12 (2022), p. 60.
\bibitem{AKRY} M. Agranovsky, A. Koldobsky, D. Ryabogin, V. Yaskin, {\em An analogue of polynomially integrable domains in even-dimensional spaces}, J. Math. Aanal. Appl., 529 (2024), 1-12.
\bibitem{ABKVY} M. Agranovsky, J. Boman, A. Koldobsky, V. Vassiliev, V. Yaskin,  {\em Algebraically integrable bodies and related properties of the Radon transform},
A. Koldobsky, A. Volberg (eds), Harmonic Analysis and Comvexity, Advances in Analysis and Geometry, 9, De Gruyter, Berlin, (2023), 1–36.
\bibitem{Ar} V. I. Arnold, {\em Arnold's Problem, 2nd edition}, Springer-Verlag, Belin, (2004).
\bibitem{ArKvant} V. I. Arnold, {\em Kepler's second law and the topology of Abelian integrals}, Kvant Selecta: Algebra and Anal. II, Math. World, vjl.15. AMS, Providence,
RI, 1999, 135-140.
\bibitem{ArVas} V.I. Arnold, V.A. Vassiliev, {\em Newton’s Principia read 300 years later}, Notices AMS, 36:9 (1989), 1148-1154.
\bibitem {B1} J. Boman, {\em A hypersurface containing the support of a Radon transform must be an ellipsoid. I: The symmetric case}, J. Geom. Anal. 31 (2021), 2726–2741;
    \bibitem {B2} J. Boman, {\em A hypersurface containing the support of a Radon transform must be an ellipsoid. II: The general case}, J. Inverse Ill-posed Probl. 29 (2021), 351–367.
%\bibitem{Vas2} V.I. Arnold, V.A. Vassiliev, {\em Addendum to [3]}, Notices AMS, 37:1, 144.
%\bibitem{Fu} W. Fulton, {\em Algebraic Topology: A First Course}, Springer Sci, 1995.
\bibitem{Gardner} R. J. Gardner, {\em Geometric Tomography}, Cambridge Univ. Press, 2006.
\bibitem{G} M. Goresky, R. MacPherson, {\em Stratified Morse Theory}, Springer Verlag, 1988.
\bibitem{He} S. Helgason,  {\em Groups and Geometric Analysis}, Acad.Press, 1984.
\bibitem{JM} J. Mather, {\em Notes on Topological Stability}, Bull. of AMS, 49 (4); 475-506 (2012).
%\bibitem{Principia} I. Newton {\em Philosophiae Naturalis Principia Mathematica}, (1687), London.
%\bibitem {IP} J.Ilmavirta, G. P. Paternain, {\em Functions of constant geodesic X-ray transform}, Inverse Problems, vol. 35, 6, 2019; %arXiv.1702.00429.
\bibitem{Principia} I. Newton, {\em Philosophiae Naturalis Principia Mathematica}, London, 1687.
\bibitem {KMY} A. Koldobsky, A. Merkurjev, V. Yaskin, {\em On polynomially integrable convex bodies}, Advances in Math., vol. 320, (2017), 876-886; arXiv.1702.00429.
\bibitem{Sp} E. H.  Spanier, {\em Algebraic Topology},McGraw-Hill book co.,  1966.
\bibitem{T} R. Thom,  {\em Ensembles et morphismes stratifiés}, Bulletin of the American Mathematical Society. 75 (2), (1969) 240–284.
\bibitem{Vas1} V. A. Vassiliev, {\em Newton's Lemma XXVIII on integrable ovals in higher dimensions and reflection groups}, Bull. Lond. Math. Soc., 47:2 (2015), 290-300.
\bibitem{Vas2} V.A. Vassiliev, {\em Applied Picard-Lefschetz Theory}, Amer. Math. Soc., Providence, RI,, (2002).
\bibitem{Vas3} V.A. Vassiliev {\em Algebroidally integrable bodies,} Arnold Math. J., 6 (2020), pp. 291-309.
\bibitem{Wiki} Wikipedia, {\em Newton's Theorem About Ovals}.
\bibitem{YM} V. Yaskin, B. Zawalski, {\em  On separably integrable symmetric convex bodies}, Adv. Math., 441 (2024), 2-23.

\end{thebibliography}
\end{document}